\documentclass[a4paper,10pt]{article}
\usepackage[latin1]{inputenc}
\usepackage[T1]{fontenc}
\usepackage{amsmath}
\usepackage{amsfonts}
\usepackage{amstext}
\usepackage{amsthm}
\usepackage[all]{xy}
\usepackage{geometry}
\usepackage{graphicx}
\usepackage{slashed}
\newtheorem{theorem}{Theorem}[section]

\newtheorem{prop}[theorem]{Proposition}
\newtheorem{cor}[theorem]{Corollary}

\DeclareMathOperator{\spin}{Spin}
\DeclareMathOperator{\End}{End}
\DeclareMathOperator{\aut}{Aut}
\DeclareMathOperator{\GL}{GL}
\DeclareMathOperator{\SO}{SO}

\title{A remark on the space of metrics having non-trivial harmonic spinors}
\author{Nils Waterstraat}

\begin{document}
\date{}
\maketitle

\begin{abstract}
Let $M$ be a closed spin manifold of dimension $n\equiv 3\mod 4$. We give a simple proof of the fact that the space of metrics on $M$ with invertible Dirac operator is either empty or it has infinitely many path components.
\end{abstract}

\section{Introduction}
Let\let\thefootnote\relax\footnotetext{2010 Mathematics Subject Classification: Primary 58J30; Secondary 53C27} $M$ be an $n$-dimensional closed spin manifold and let $R(M)$ be the space of all Riemannian metrics on $M$. For any choice of a metric $g\in R(M)$, we can build the associated spinor bundle $\Sigma_gM$ and obtain a natural first order operator $D_g$ acting on sections of $\Sigma_gM$ and which we call the \textit{Dirac operator}. Elements of $\ker D_g$ are called \textit{harmonic spinors} and their existence has been studied for a long time. While one can show that on $S^2$ no non-trivial harmonic spinors exist for any choice of $g$ (cf. \cite{BaerEstimates}), it is conjectured that on every closed spin manifold of dimension $n\geq 3$ there exists a Riemannian metric $g$ such that $\ker D_g\neq 0$. The conjecture has been proved by N. Hitchin in \cite{Hitchin} if $n\equiv 0,\pm 1\mod 8$ and by C. Bär in \cite{BaerHarmonic} if $n\equiv 3\mod 4$.\\
As a more general question, one may ask how many metrics exist on $M$ such that the corresponding Dirac operator has non-trivial kernel. A possible way to study this question is to consider the complementary set of metrics $R^\textup{inv}(M)$ consisting of all metrics $g\in R(M)$ such that $\ker D_g=0$. M. Dahl showed in \cite{Dahl} that elements of $R^\textup{inv}(M)$ can be extended to $R^\textup{inv}(W)$ if $W$ is the trace of a surgery of codimension at least $3$ on $M$. By using the Atiyah-Singer index theorem and special metrics on the spheres originating from the study of positive scalar curvature, he concluded from this result that $R^\textup{inv}(M)$ is in  all dimensions $n\geq 5$ which were considered by Hitchin and B\"ar either empty or disconnected. Moreover, in the case $n\equiv 3\mod 4$, $n\geq 7$, he even obtained that, if non-empty, $R^\textup{inv}(M)$ has infinitely many path components. Recently he improved this conclusion in collaboration with N. Grosse to dimension 3 by studying extensions of metrics to attached handles (cf. \cite{DahlGrosse}).\\
The aim of this article is to show that the existence of infinitely many connected components of $R^\textup{inv}(M)$ in all dimensions $n\equiv 3\mod 4$ can be derived easily from B\"ar's results in \cite{BaerHarmonic} by using spectral flow and rather elementary homotopy arguments.\\
Finally, we want to mention that B\"ar improved his theorem in \cite{BaerTwisted} to twisted Dirac operators. Note that for any fixed pair $(F,\nabla)$ of a bundle $F$ over $M$ and a connection $\nabla$ on $F$, we obtain a family $D^{(F,\nabla)}$ of twisted Dirac operators which is again parametrised by the space of Riemannian metrics $R(M)$ on $M$. We believe that one can extend our argument here to this case by using the results from \cite{BaerTwisted} instead of \cite{BaerHarmonic}. Accordingly, we conjecture that the corresponding space $R^\textup{inv}_{(F,\nabla)}(M)=\{g\in R(M):\ker D^{(F,\nabla)}_g=0\}$ is either empty or it has infinitely many path components.\\
\textbf{Acknowledgement}: I want to thank Jacobo Pejsachowicz for encouraging me to write this paper and for several helpful suggestions improving its presentation.


\section{Preliminaries: Dirac operators}
In this section we recall briefly the definition of spinor bundles and their Dirac operators. Among the many references for these topics we want to mention \cite{Hijazi} and \cite{AmmannTori}, on which we base our exposition. In order to simplify the presentation we assume throughout that $M$ is an oriented closed manifold of odd dimension $n\geq 3$.\\
We denote by $\GL^+(M)$ the principal $\GL^+(n;\mathbb{R})$-bundle of oriented bases over $M$ and recall that $\GL^+(n;\mathbb{R})$ has a unique connected 2-fold covering $\Theta:\widetilde{\GL^+}(n;\mathbb{R})\rightarrow\GL^+(n;\mathbb{R})$ since the fundamental group of $\GL^+(n;\mathbb{R})$ is of order two. A \textit{spin structure} on $M$ is a pair $(\widetilde{\GL^+}(M),\vartheta)$, where $\widetilde{\GL^+}(M)$ is a principal $\widetilde{\GL^+}(n;\mathbb{R})$-bundle over $M$ and $\vartheta:\widetilde{\GL^+}(M)\rightarrow\GL^+(M)$ is a 2-fold covering such that

\begin{align*}
\pi\circ\vartheta=\widetilde{\pi}\quad\text{and}\quad \vartheta(u\cdot v)=\vartheta(u)\cdot\Theta(v),\quad\text{for all}\,\, v\in\widetilde{\GL^+}(n;\mathbb{R}),\, u\in\widetilde{\GL^+}(M),
\end{align*}   
where $\pi$ and $\widetilde{\pi}$ denote the corresponding projections of the bundles. Henceforth we assume that $M$ is a \textit{spin manifold}, that is, $M$ is oriented and a spin structure on $M$ is given. Note that so far we have not required that $M$ is endowed with a Riemannian metric.\\
Let now $g$ be a Riemannian metric on $M$ and denote by $\SO(M,g)$ the associated principal $\SO(n)$-bundle of positively oriented orthonormal bases. Then $\spin(M,g):=\vartheta^{-1}(\SO(M,g))$ is a principal $\spin(n)$-bundle over $M$, where $\spin(n):=\Theta^{-1}(\SO(n))$ is the unique connected 2-fold covering of $\SO(n)$. Let $\rho:\mathbb{C}l_n\rightarrow\End(\Sigma_n)$ denote the usual irreducible representation of the complex Clifford algebra, where $\Sigma_n$ is the space of complex spinors. We fix an inner product $\langle\cdot,\cdot\rangle$ on $\Sigma_n$ such that $\langle\rho(x)\sigma_1,\rho(x)\sigma_2\rangle=\langle\sigma_1,\sigma_2\rangle$ for all $x\in\mathbb{R}^n$, $\|x\|=1$, and $\sigma_1,\sigma_2\in\Sigma_n$. If now $\rho':\spin(n)\rightarrow\aut(\Sigma_n)$ denotes the complex spinor representation of $\spin(n)$, which is obtained by restricting $\rho$ to $\spin(n)\subset\mathbb{C}l_n$, then the \textit{spinor bundle} $\Sigma_g M$ of $M$ with respect to $g$ is defined as the associated vector bundle $\spin(M,g)\times_{\rho'}\Sigma_n$.\\
The representation $\rho$ induces a \textit{Clifford multiplication} on $\Sigma_g M$, that is, a complex linear vector bundle homomorphism

\begin{align*}
m:T^\ast M\otimes\Sigma_g M\rightarrow\Sigma_g M,\quad X^\flat\otimes\varphi\mapsto X\cdot\varphi
\end{align*}      
such that $X\cdot(Y\cdot\varphi)+Y\cdot(X\cdot\varphi)=-2g(X,Y)\varphi$ for all $X,Y\in TM$ and $\varphi\in\Sigma_g M$. Moreover, the inner product on $\Sigma_n$ gives rise to an Hermitian structure on the bundle $\Sigma_g M$ such that $\langle X\cdot\varphi,\psi\rangle=-\langle\varphi,X\cdot\psi\rangle$ for all $X\in TM$ and $\varphi,\psi\in\Sigma_g M$. Finally, the Levi-Civita connection on $TM$ induces a connection on $\SO(M,g)$ and this connection lifts in a canonical way to a connection on $\spin(M,g)$. The associated covariant derivative $\nabla:C^\infty(M,\Sigma_g M)\rightarrow C^\infty(M,T^\ast M\otimes\Sigma_g M)$ on the spinor bundle has the properties

\begin{align*}
X\langle\varphi,\psi\rangle=\langle\nabla_X\varphi,\psi\rangle+\langle\varphi,\nabla_X\psi\rangle\quad \text{and}\quad \nabla_X(Y\cdot\varphi)=(\nabla^{TM}_XY)\cdot\varphi+Y\cdot(\nabla_X\varphi)
\end{align*}
for vector fields $X,Y$ and a spinor field $\varphi$.\\
Now the Dirac operator with respect to the metric $g$ is defined by 

\begin{align*}
D_g=m\circ\nabla:C^\infty(M,\Sigma_g M)\rightarrow C^\infty(M,\Sigma_g M)
\end{align*}
and is an elliptic, essentially selfadjoint differential operator of first order.


\section{The Proof}
We assume from now on that $M$ is a closed spin manifold of dimension $n\equiv 3\mod 4$. We denote by $R(M)$ the space of all Riemannian metrics on $M$ with the $C^1$-topology and note that it is obviously contractible. Moreover, we define

\begin{align*}
R^\textup{inv}(M)=\{g\in R(M):\,\ker D_g=0\}\subset R(M)
\end{align*} 
and recall that our aim is to show that this set has infinitely many path components if it is not empty. Accordingly, we assume henceforth that $R^\textup{inv}(M)\neq\emptyset$ and now we conclude in three steps the announced disconnectedness of this space.


\subsection*{Step 1: The spectral flow}
Since our operators $D_g$, $g\in R(M)$, are essentially selfadjoint, they have real spectra. Moreover, by ellipticity their spectra are discrete and consist entirely of eigenvalues of finite multiplicity. We define for any compact interval $[a,b]\subset\mathbb{R}$ a non-negative integer by

\begin{align*}
m(g,[a,b])=\sum_{\lambda\in [a,b]}\dim\ker(D_g-\lambda\cdot id).
\end{align*}
Next we quote the following stability result for the spectra of the operators $D_g$ that can be found in \cite[Prop. 7.1]{BaerHarmonic}.

\begin{theorem}
Let $(M,g)$ be a closed Riemannian spin manifold with Dirac operator $D_g$. Let $\varepsilon>0$ and let $\Lambda>0$ such that $-\Lambda,\Lambda\notin\sigma(D_g)$. Write

\begin{align*}
\sigma(D_g)\cap(-\Lambda,\Lambda)=\{\lambda_1\leq\lambda_2\leq\ldots\leq\lambda_k\}.
\end{align*}
Then there exists a neighbourhood of $g$ in the $C^1$-topology such that for any metric $\tilde{g}$ in this neighbourhood with Dirac operator $D_{\tilde{g}}$ the following holds:

\begin{itemize}
\item $\sigma(D_{\tilde{g}})\cap(-\Lambda,\Lambda)=\{\mu_1\leq\mu_2\leq\ldots\leq\mu_k\}$,
\item $|\lambda_i-\mu_i|<\varepsilon$, $i=1,\ldots,k$.
\end{itemize}
The eigenvalues $\lambda_i$ and $\mu_i$ are repeated according to their multiplicities.
\end{theorem}

We obtain immediately the following corollary.

\begin{cor}\label{cor-stabspec}
For all $g_0\in R(M)$ and $\Lambda>0$ such that $\pm\Lambda\notin\sigma(D_{g_0})$ there exists an open neighbourhood $N(g_0,\Lambda)\subset R(M)$ such that $\pm\Lambda\notin\sigma(D_g)$ and $m(g,[-\Lambda,\Lambda])=m(g_0,[-\Lambda,\Lambda])$ for all $g\in N(g_0,\Lambda)$.
\end{cor}

Let now $\gamma:I\rightarrow R(M)$ be a path of metrics. Because of corollary \ref{cor-stabspec} we can find a decomposition $0=t_0< t_1<\ldots<t_N=1$ and positive numbers $a_1,\ldots,a_N$ such that the functions $[t_{i-1},t_i]\ni t\mapsto m(\gamma(t),[-a_i,a_i])$ are constant. We define

\begin{align}\label{sfl}
\Gamma(\gamma)=\sum^N_{i=1}{m(\gamma(t_i),[0,a_i])-m(\gamma(t_{i-1}),[0,a_i])}\in\mathbb{Z}
\end{align}  
and note that, roughly speaking, $\Gamma(\gamma)$ counts the number of negative eigenvalues of $D_{\gamma(0)}$ that become positive as the parameter $t$ travels from $0$ to $1$ minus the number of positive eigenvalues of $D_{\gamma(0)}$ that become negative; i.e., the net number of eigenvalues which cross zero. The formula \eqref{sfl} corresponds precisely to the definition of the spectral flow for paths of selfadjoint Fredholm operators acting on a fixed Hilbert space which can be found for example in \cite{Phillips} and \cite{UnbSpecFlow}. Accordingly, one can show verbatim as in \cite{Phillips} that $\Gamma(\gamma)$ indeed does only depend on the path $\gamma$ and not on the choices of the $t_i,a_i$, $i=1,\ldots N$. Moreover, if $\gamma,\tilde{\gamma}:I\rightarrow R(M)$ are two paths of metrics, then the following properties hold:

\begin{itemize}
	\item[i)] $\Gamma(\gamma)=0$ if $\gamma(t)\in R^\textup{inv}(M)$ for all $t\in[0,1]$,
	\item[ii)] $\Gamma(\gamma\ast\tilde{\gamma})=\Gamma(\gamma)+\Gamma(\tilde{\gamma})$, whenever the concatenation $\gamma\ast\tilde{\gamma}$ exists,
	\item[iii)] $\Gamma(\gamma^{-1})=-\Gamma(\gamma)$, where $\gamma^{-1}(t)=\gamma(1-t)$, $t\in I$,
	\item[iv)] $\Gamma(\gamma)=\Gamma(\tilde{\gamma})$ if $\gamma\simeq\tilde{\gamma}$ through a homotopy having ends in $R^\textup{inv}(M)$. 
\end{itemize}

Note that the first three properties are immediate consequences of the definition. The homotopy invariance can be obtained again verbatim as in \cite{Phillips}.


\subsection*{Step 2: The range of $\Gamma$}
Our argument in this section is based on results from \cite{BaerHarmonic} which we introduce before we proceed with the proof. At first, we need the existence of the following metrics on the sphere $S^n$, that were constructed in \cite[\S3]{BaerHarmonic}.  

\begin{prop}\label{prop-existence}
For $n\equiv 3\mod 4$ and any integer $m>0$, there exists a path of metrics $g^{m}_t$, $t\in[0,1]$, on $S^n$ such that the following holds for the associated Dirac operators $\slashed{D}^m_t$:

\begin{itemize}
\item there is $\lambda(t)\in\sigma(\slashed{D}^m_t)$ such that $\lambda(0)=-1$ and $\lambda(1)=1$,
\item $\lambda(t)$ depends linearly on $t$,
\item the multiplicity of $\lambda(t)$ is constant in $t$ and greater than $m$,
\item $\lambda(t)$ is the only eigenvalue of $\slashed{D}^m_t$ in the interval $[-2,2]$.
\end{itemize}
\end{prop}

B\"ar combined in \cite{BaerHarmonic} proposition \ref{prop-existence} and a general gluing theorem for Dirac operators \cite[theorem B]{BaerHarmonic} to conclude the existence of non-trivial harmonic spinors in dimensions $n\equiv 3\mod 4$. Actually, in order to find the spinors he just needed a special case of his gluing theorem which reads as follows. 

\begin{theorem}\label{theorem-baer}
Let $(M,g)$ be a closed Riemannian spin manifold of odd dimension $n\geq 3$. Let $D_g$ be the corresponding Dirac operator and let $\slashed{D}$ denote the Dirac operator on $S^n$ with respect to some Riemannian metric. Finally, let $\Lambda>0$ be such that $\pm\Lambda\notin\sigma(D_g)\cup\sigma(\slashed{D})$. Write

\begin{align*}
(\sigma(D_g)\cup\sigma(\slashed{D}))\cap(-\Lambda,\Lambda)=\{\lambda_1\leq\lambda_2\leq\ldots\leq\lambda_k\}.
\end{align*}
Then for any $\varepsilon>0$ there exists a Riemannian metric $\tilde{g}$ on $M$ such that the corresponding Dirac operator $D_{\tilde{g}}$ has the following properties: 

\begin{itemize}
	\item[i)] $\pm\Lambda\notin\sigma(D_{\tilde{g}})$,
	\item[ii)] $\sigma(D_{\tilde{g}})\cap(-\Lambda,\Lambda)=\{\mu_1\leq\mu_2\leq\ldots\leq\mu_k\}$
	\item[iii)] $|\lambda_j-\mu_j|<\varepsilon$ for $j=1,\ldots,k$.
\end{itemize}
The eigenvalues $\lambda_i$ and $\mu_i$ are repeated according to their multiplicities.
\end{theorem}

We now take some metric $g_0\in R^\textup{inv}(M)$. Because of the conformal covariance of the Dirac operator (cf. \cite[Prop. 5.13]{Hijazi}), we can assume that  $[-2,2]\cap\sigma(D_{g_0})=\emptyset$ simply by rescaling the metric if necessary.\\
Let $m>0$ be an integer and consider the operators $\slashed{D}^m_t$ on $S^n$ from proposition \ref{prop-existence}. Recall that we denote by $\lambda(t)$ the unique eigenvalue of $\slashed{D}^m_t$ in the interval $[-2,2]$ and that $\lambda(t)$ depends linearly on $t$ with $\lambda(0)=-1$, $\lambda(1)=1$.\\
We now apply theorem \ref{theorem-baer} for $\Lambda=2$ and $\varepsilon=\frac{1}{2}$ to $D_{g_0}$ and the operators $\slashed{D}^m_t$, $t\in[0,1]$, on $S^n$. We obtain for any $t\in[0,1]$ a metric $\tilde{g}_t$ on $M$ such that each eigenvalue of $D_{\tilde{g}_t}$ in the interval $[-2,2]$ is of distance less then $\frac{1}{2}$ to $\lambda(t)$. In particular, $D_{\tilde{g}_0}$ and $D_{\tilde{g}_1}$ are invertible and hence $\{\tilde{g}_t\}_{t\in[0,1]}$ defines a path $\gamma:(I,\partial I)\rightarrow(R(M),R^\textup{inv}(M))$. Moreover, the function $t\mapsto m(\gamma(t),[-2,2])$ is constant on the whole interval $[0,1]$. Hence we finally obtain from the definition of $\Gamma$

\begin{align*}
\Gamma(\gamma)&=m(\tilde{g}_1,[0,2])-m(\tilde{g}_0,[0,2])=m(\tilde{g}_1,[0,2])=\dim\ker(\slashed{D}^m_1-id)>m.
\end{align*}
To sum up, we have shown that the set

\begin{align*}
\{\Gamma(\gamma):\,\,\gamma:(I,\partial I)\rightarrow(R(M),R^\textup{inv}(M))\,\text{continuous}\}\subset\mathbb{Z}
\end{align*}
is not bounded from above.


\subsection*{Step 3: The final argument}
We fix some $g_0\in R^\textup{inv}(M)$. Our first aim of this final step is to construct inductively a sequence of paths $\gamma_k:(I,\partial I)\rightarrow(R(M),R^\textup{inv}(M))$, $k\in\mathbb{N}$, such that $\gamma_k(0)=g_0$ for all $k\in\mathbb{N}$ and $\Gamma(\gamma_i)\neq\Gamma(\gamma_j)$ for all $i\neq j$.\\
Let $\gamma_1$ be the constant path $\gamma_1\equiv g_0\in R^\textup{inv}(M)$. Assume that we have already constructed $\gamma_1,\ldots,\gamma_k:(I,\partial I)\rightarrow (R(M),R^\textup{inv}(M))$ such that $\gamma_i(0)=g_0$, $i=1,\ldots k$, and $\Gamma(\gamma_i)\neq\Gamma(\gamma_j)$ for all $i\neq j$.\\
According to the second step of our proof we can find a path $\tilde{\gamma}:(I,\partial I)\rightarrow(R(M),R^\textup{inv}(M))$ such that 

\begin{align}\label{formulaend}
\Gamma(\tilde{\gamma})>\max_{1\leq i,j\leq k}|\Gamma(\gamma_i)-\Gamma(\gamma_j)|.
\end{align}
Moreover, we choose a path $\hat{\gamma}:(I,\partial I)\rightarrow(R(M),R^\textup{inv}(M))$ such that $\hat{\gamma}(0)=g_0$ and $\hat{\gamma}(1)=\tilde{\gamma}(0)$. Then $ \hat{\gamma}\ast\tilde{\gamma}:(I,\partial I)\rightarrow(R(M),R^\textup{inv}(M))$ and we set $\gamma_{k+1}=\hat{\gamma}\ast\tilde{\gamma}$ if $\Gamma(\hat{\gamma}\ast\tilde{\gamma})\neq \Gamma(\gamma_j)$ for all $j=1,\ldots,k$.\\
If, on the other hand, $\Gamma(\hat{\gamma}\ast\tilde{\gamma})= \Gamma(\gamma_j)$ for some $j=1,\ldots,k$, then we set $\gamma_{k+1}=\hat{\gamma}$. In order to justify this choice, assume that also $\Gamma(\hat{\gamma})=\Gamma(\gamma_i)$ for some $1\leq i\leq k$. Then we obtain

\begin{align*}
\Gamma(\gamma_j)&=\Gamma(\hat{\gamma}\ast\tilde{\gamma})=\Gamma(\hat{\gamma})+\Gamma(\tilde{\gamma})=\Gamma(\gamma_i)+\Gamma(\tilde{\gamma}),
\end{align*}  
which contradicts \eqref{formulaend}. Hence we indeed obtain a sequence $\{\gamma_k\}_{k\in\mathbb{N}}$ with the required properties.\\
We now finish our proof by claiming that the metrics $\gamma_k(1)$, $k\in\mathbb{N}$, all lie in different path components of $R^\textup{inv}(M)$. Assume on the contrary that we can find $i,j\in\mathbb{N}$, $i\neq j$, and a path $\tilde{\gamma}:I\rightarrow R^\textup{inv}(M)$ such that $\tilde{\gamma}(0)=\gamma_i(1)$ and $\tilde{\gamma}(1)=\gamma_j(1)$. Then $\gamma_i\ast\tilde{\gamma}\ast\gamma^{-1}_j$ is a closed path with initial point $g_0\in R^\textup{inv}(M)$. Since $R(M)$ is contractible, $\gamma_i\ast\tilde{\gamma}\ast\gamma^{-1}_j$ is homotopic to the constant path $\gamma_1\equiv g_0$ through a $g_0$-preserving homotopy. We obtain from the properties of $\Gamma$

\begin{align*}
0=\Gamma(\gamma_1)=\Gamma(\gamma_i\ast\tilde{\gamma}\ast\gamma^{-1}_j)=\Gamma(\gamma_i)+\Gamma(\tilde{\gamma})+\Gamma(\gamma^{-1}_j)=\Gamma(\gamma_i)+\Gamma(\gamma^{-1}_j)
\end{align*}
and hence $\Gamma(\gamma_i)=\Gamma(\gamma_j)$ contradicting the construction of the sequence $\{\gamma_k\}_{k\in\mathbb{N}}$.

\vspace{1cm}

Dipartimento di Scienze Matematiche\\
Politecnico di Torino\\
Corso Duca degli Abruzzi, 24\\
10129 Torino, Italy\\
E-mail: waterstraat@daad-alumni.de

\end{document}